\documentclass[a4paper,notitlepage, reqno]{amsart}
\usepackage{fullpage}
\usepackage[colorlinks=true,linkcolor=blue,unicode, psdextra]{hyperref}
\usepackage{color}
\usepackage{amsthm}
\usepackage{amsmath}
\usepackage{amssymb}
\usepackage{enumitem}
\usepackage{mathtools}
\usepackage{mathrsfs}
\usepackage[utf8]{inputenc}

\newtheorem*{theorem*}{Theorem}
\newtheorem*{proposition*}{Proposition}
\newtheorem*{lemma*}{Lemma}
\newtheorem*{corollary*}{Corollary}
\newtheorem*{remark*}{Remark}

\newtheorem{theorem}{Theorem}

\theoremstyle{definition}

\newtheorem{remark}[theorem]{Remark}

\allowdisplaybreaks

\date{}
\title{On the origin of the Jacobian conjecture}
\author{Lázaro Orlando Rodríguez Díaz}
\address{Instituto de Matem\'{a}tica, 
	Universidade Federal do Rio de Janeiro, RJ, Brazil.}
\email{lazarord@im.ufrj.br}	

\begin{document}
	
\begin{abstract}
	The Jacobian conjecture is thought to have been proposed by O. H. Keller in 1939. However, we have found that the statement of the conjecture is precisely the main result of a paper published by L. Kraus in 1884. Although the final step of Kraus's proof is flawed, the ideas he introduced anticipated approaches to the problem that would only emerge more than a century later. Interestingly, the root of Kraus's error remains the principal obstacle to algebro-geometric approaches: controlling the ramification at infinity.   
\end{abstract}
\maketitle

\section*{}
It is believed that the Jacobian conjecture was first proposed by O. H. Keller \cite{Keller39} in 1939; for a comprehensive historical review, see Bass, Connell, and Wright \cite{Bass_Connell_Wright_82}. However, a recent search of the zbMATH Open database for the German term \textit{Functionaldeterminanten} led to the discovery of a paper by L. Kraus \cite{Kraus1884} dating back to 1884.  While the reviewer’s summary \footnote{\url{https://zbmath.org/16.0232.03}} alone suggests that Kraus may have addressed the conjecture, a close study of the article \footnote{\url{https://archive.org/details/sitzungsbericht398klasgoog/page/821/}} reveals that he sought to prove the exact statement now known as the plane Jacobian conjecture. Quoting Kraus himself:
\begin{quote}
``The final result is therefore as follows:

If two polynomials of the variables $x, y$
\begin{align*}
F(xy), G(xy)
\end{align*}
have the property that their functional determinant
\begin{align*}
\Delta(xy)=\begin{vmatrix}
\frac{\partial F(xy)}{\partial x} & \frac{\partial F(xy)}{\partial y}\\
\frac{\partial G(xy)}{\partial x} & \frac{\partial G(xy)}{\partial y}
\end{vmatrix}
\end{align*}
is a constant different from zero, then from the system of equations
\begin{align*}
F(xy)-\eta=0\\
G(xy)-\xi=0
\end{align*}
the quantities $x, y$ are obtained as polynomials of the variables $\xi, \eta$.

And conversely: if from that system of equations the quantities $x, y$ are obtained as polynomials of the variables $\xi, \eta$, then $\Delta(xy)$ is a constant different from zero." \footnote{
	``Das Schlussresultat ist also folgendes:
	
	Haben zwei ganze, rationale Functionen der Grössen $x, y$
	\begin{align*}
	F(xy), G(xy)
	\end{align*}
	die Eigenschaft, dass ihre Functionaldeterminante
	\begin{align*}
	\Delta(xy)=\begin{vmatrix}
	\frac{\partial F(xy)}{\partial x} & \frac{\partial F(xy)}{\partial y}\\
	\frac{\partial G(xy)}{\partial x} & \frac{\partial G(xy)}{\partial y}
	\end{vmatrix}
	\end{align*}
	eine von Null verschiedene Constante ist, so ergeben sich aus dem Gleichungssysteme
	\begin{align*}
	F(xy)-\eta=0\\
	G(xy)-\xi=0
	\end{align*}
	die Grössen $x, y$ als ganze rationale Functionen der Variablen $\xi, \eta$.
	
	Und umgekehrt: Ergeben sich aus jenem Gleichungssystem die Grössen $x, y$ als ganze, rationale Functionen der Variablen $\xi, \eta$, so ist $\Delta(xy)$ eine von Null verschiedene Constante."} \cite[p. 825]{Kraus1884} 	
\end{quote}

Kraus's approach is based on: the theory of the resultant and the theory of algebraic functions in complex analysis. The use of the resultant in this context is well known to experts (cf. Adjamagbo and van den Essen \cite{Adjamagbo_Essen90}, Chadzyński and Krasi{\'n}ski \cite{Chadzynski_Krasinski_92}, McKay and Wang \cite{McKay_Wang86}, Moh \cite[Section 7]{Moh98}). The novelty of Kraus’s strategy, however, lies in his use of the irreducibility of the fibers (or the irreducibility of the resultant). This perspective allowed him to take advantage of the theory of algebraic functions. Notably, modern results by Kaliman \cite{Kaliman93} confirm that one may assume the irreducibility of the fibers of a coordinate polynomial without loss of generality when addressing the conjecture.

Although the final step of Kraus’s proof is flawed, his methods anticipated modern efforts by over a century. Remarkably, the core of his error remains the primary obstacle in algebro-geometric approaches: controlling the ramification at infinity.

Kraus's work \cite{Kraus1884} was preceded by two papers \footnote{\url{https://archive.org/details/sitzungsberichte1882knig/page/338/}} \footnote{\url{https://archive.org/details/sitzungsberichte1883knig/page/187/}} exploring the invertibility of polynomial maps \cite{Kraus1882, Kraus1883}. It is worth noting that Kraus's title \cite{Kraus1884} mirrors Jacobi's 1841 work \cite{Jacobi1841, Jacobi1896}, in which Jacobi detailed the properties of what is now termed the \emph{Jacobian determinant}. While Kraus's paper is cited in Muir’s historical treatise \cite[p. 254]{Muir1923}, it has remained overlooked in the context of the conjecture. For a brief biographical account of Ludwig Kraus (1857-1886), see \cite{Weyr1886}.

We emphasize that Kraus's paper \cite{Kraus1884} is notably clear and detailed, especially considering its nineteenth-century origin. While certain steps in his argument require further justification or do not go through as stated without modification, we do not attempt here to rectify Kraus's proof in full detail. Our aim is instead to reconstruct the argument using results that are now available and to pinpoint precisely where the critical flaw arises. 

\section{Kraus's proof}

In this section we reconstruct Kraus's proof. We essentially follow Kraus's strategy, while using modern results to justify several claims. Along the way, we will see how Kraus’s ideas resurfaced more than a century later in subsequent attempts to address this problem.

Let $(p,q): \mathbb{C}^{2}\to \mathbb{C}^{2}$ be a polynomial map whose Jacobian $J(p,q)=\frac{\partial p}{\partial x}\frac{\partial q}{\partial y}-\frac{\partial p}{\partial y}\frac{\partial q}{\partial x}$ is equal to $1$. We aim to prove that this map is invertible. To this end, we are guided by a simple resultant criterion for invertibility of polynomial maps due to Adjamagbo and van den Essen. Let $u$ and $v$ be indeterminates over the polynomial ring $\mathbb{C}[x,y]$.

\begin{theorem}[\cite{Adjamagbo_Essen90}]\label{thm-Adjamagbo_Essen}
	Let $(f_1(x,y),f_2(x,y))$ be a polynomial map of $\mathbb{C}^{2}$. There is an equivalence between:
	\begin{itemize}
		\item[i)] $(f_1,f_2)$ is invertible.
		\item[ii)] There exist $\lambda_{1}, \lambda_{2}\in \mathbb{C}^{*}$ and $g_{1}, g_{2}\in \mathbb{C}[u,v]$ such that $$\operatorname{Res}_{y}(f_{1}-u, f_{2}-v)=\lambda_{1}(x-g_{1}) \quad \text{and} \quad \operatorname{Res}_{x}(f_{1}-u, f_{2}-v)=\lambda_{2}(y-g_{2}).$$
	\end{itemize}
	Furthermore, if $(f_1,f_2)$ is invertible, then $(g_1,g_2)$ is its inverse.
\end{theorem}

To prove that $(p,q)$ is invertible, we need less than what was required in the previous theorem, due to the extra hypothesis that $J(p,q)=1$. By composing with suitable invertible linear transformations, we can always assume that:
\begin{equation}
\begin{aligned}\label{eq_form_of_p_and_q}
p(x,y)=y^{\operatorname{deg}p}+\text{lower terms in}\;y, \\ q(x,y)=y^{\operatorname{deg}q}+\text{lower terms in}\;y.
\end{aligned}
\end{equation}
 
In order to prove that $(p,q)$ is invertible, we just need to consider the polynomial
\begin{align}\label{eq_Resultant}
\operatorname{Res}_{y}(p-u, q-v)=r_{n}(u,v)x^{n}+\cdots+r_{0}(u,v),
\end{align}
and show that $n=1$. Because of (\ref{eq_form_of_p_and_q}) and $J(p,q)=1$, we always have $n\ge 1$ and $r_{0}(u,v)\ne 0$, see Sakkalis \cite[Theorem 1]{Sakkalis93}. Suppose we are able to prove that $n=1$; then, specializing $u\to p$ and $v\to q$ in (\ref{eq_Resultant}), we obtain $r_{1}(p,q)x+r_{0}(p,q)=0$, that is, $x\in\mathbb{C}(p,q)$. On the other hand, by a theorem of Formanek \cite[Theorem 2]{Formanek94}, we know that a polynomial map $(p,q)$ with constant Jacobian satisfies $\mathbb{C}(x,y)=\mathbb{C}(x,p,q)$. Combining $x\in\mathbb{C}(p,q)$ and $\mathbb{C}(x,y)=\mathbb{C}(x,p,q)$, we obtain $\mathbb{C}(x,y)=\mathbb{C}(p,q)$; that is, the map $(p,q)$ is birational. That $n=1$ implies $(p,q)$ is birational also follows from the fact that, under the condition (\ref{eq_form_of_p_and_q}), $n$ equals the geometric degree of the dominant map $(p,q)$; see P\l oski \cite[Lemma 3.1]{Ploski94}. However, birational polynomial maps with constant Jacobian are invertible, as proved by Keller \cite[Satz 3]{Keller39}. 

\begin{remark}\label{rmk_specialization}
	In general, it is not true that the specialization of the resultant is equal to the resultant of the specialized polynomials. In our case, since we are assuming that the map $(p,q)$ has the form given in (\ref{eq_form_of_p_and_q}), we are allowed to specialize, as proved in \cite[Chapter 3, Section 6, Proposition 6]{Cox_Little_Shea_2025}.
\end{remark}

We now turn to the examination of the polynomial $\operatorname{Res}_{y}(p-u, q-v)$. In this context, we have the following interesting result by Abhyankar \cite[pp. 67-73]{Abhyankar71} and McKay and Wang \cite[Theorem 1]{McKay_Wang86}. For a proof of the more general case where $f(t), g(t)\in k(t)$ see Gutierrez, Rubio, and Yu \cite[Theorem 1.3]{Gutierrez_Rubio_Yu02}.  

\begin{theorem}[\cite{Abhyankar71}, \cite{McKay_Wang86}]\label{thm_Abhyankar_McKay_Wang}
	Let $k$ be a field, $t, u, v$ be indeterminates over $k$ and $f(t), g(t)\in k[t]$ not both constant. Then there exist an irreducible polynomial $m(u,v)\in k[u,v]$ and $\alpha\in k^{*}$ such that:
	\begin{align*}
	\operatorname{Res}_{t}(f(t)-u, g(t)-v)=\alpha\left(m(u,v)\right)^{l},
	\end{align*}
	where $l=[k(t):k(f,g)]$, the degree of field extensions.
\end{theorem}

To apply the above theorem to $\operatorname{Res}_{y}(p-u, q-v)$, we set $k=\mathbb{C}(x)$, and regard $p$ and $q$ as polynomials in $y$ with coefficients in $\mathbb{C}[x]$. By the previously mentioned result of Formanek \cite[Theorem 2]{Formanek94}, we have that $\mathbb{C}(x,y)=\mathbb{C}(x,p,q)$; thus, $l=[\mathbb{C}(x,y):\mathbb{C}(x,p,q)]=1$. It follows that $\operatorname{Res}_{y}(p-u, q-v)=G(x)H(x,u,v)$, where $H(x,u,v)\in \mathbb{C}[x,u,v]$ is irreducible and $G(x)\in \mathbb{C}[x]$. Moreover, it is known that $G(x)$ is constant if and only if $p$ and $q$ are monic as polynomials in $y$, see Formanek \cite[p. 375]{Formanek94} and Cassou-Noguès and Miyanishi \cite[Theorem 1.1]{Cassou-Nogues_Miyanishi09}. We conclude from (\ref{eq_form_of_p_and_q}) that
\begin{align}\label{eq-general_resultant}
\operatorname{Res}_{y}(p-u, q-v)=R(x,u,v),
\end{align}
where $R(x,u,v)\in\mathbb{C}[x,u,v]$ is an irreducible polynomial.

Kaliman proved in \cite{Kaliman93}  that to show that $(p,q)$ is invertible we can always assume that for every $c\in \mathbb{C}$ the fiber $\{(x,y)| p(x,y)=c\}$ is irreducible. Note that we can assume $p$ and $q$ take the form (\ref{eq_form_of_p_and_q}) without affecting the assumption that the fibers $\{(x,y) \mid p(x,y)=c\}$ are irreducible for every $c\in \mathbb{C}$.

In light of this, the following result from Hejmej \cite[Theorem 2.1]{Hejmej18} and García-Barroso and Gwo\'{z}dziewicz \cite[Theorem 2.5]{Garcia_Barroso-Gwozdziewic22} is noteworthy. Kraus noted \cite[p. 819]{Kraus1884} that this result was well known. In fact, we found the theorem was proved by Weierstrass \cite[pp. 47-48]{Weierstrass_1902} \footnote{\url{https://archive.org/details/mathematischewer04weieuoft/page/48/}}.

\begin{theorem}[\cite{Weierstrass_1902}, \cite{Hejmej18}, \cite{Garcia_Barroso-Gwozdziewic22}]\label{thm-Hejmej_Garcia_Barroso-Gwozdziewic}
	Let $k$ be a field, $t, v$ be indeterminates over $k$ and $f(t), g(t)\in k[t]$ be monic polynomials.  If $f(t)$ is irreducible in $k[t]$ then $\operatorname{Res}_{t}(f(t), g(t)-v)$ is a power of an irreducible polynomial in $k[v]$.
\end{theorem}

We can determine the value of the above power, and it is not surprising that it is exactly the degree of the field extension $[k(t):k(f,g)]$, as in Theorem \ref{thm_Abhyankar_McKay_Wang}. In fact, by reviewing the proof of Theorem \ref{thm_Abhyankar_McKay_Wang}, we observe that the role of the variable $u$ in $f(t)-u$ is to ensure that the polynomial $f(t)-u$ is irreducible in $k[u][t]$. The role formerly played by the variable $u$ is now replaced by the assumption in Theorem \ref{thm-Hejmej_Garcia_Barroso-Gwozdziewic} that $f(t)$ is irreducible in $k[t]$ from the outset. Although both proofs are based on field theory and the core ideas are similar, Theorem \ref{thm-Hejmej_Garcia_Barroso-Gwozdziewic} does not follow from Theorem \ref{thm_Abhyankar_McKay_Wang} by a simple specialization of the resultant.

In order to use the aforementioned theorem in our context, we let $k=\mathbb{C}(x)$, and regard $p$ and $q$ as polynomials in $y$ with coefficients in $\mathbb{C}[x]$. We know that the power equals $[\mathbb{C}(x,y):\mathbb{C}(x,p,q)]=1$. Since $p(x,y)-c=0$ is irreducible for every $c\in \mathbb{C}$, it follows from Theorem \ref{thm-Hejmej_Garcia_Barroso-Gwozdziewic} that $\operatorname{Res}_{y}(p-c, q-v)$ is an irreducible polynomial in $\mathbb{C}(x)[v]$ for all $c\in \mathbb{C}$. Furthermore, we already know from (\ref{eq-general_resultant}) that $\operatorname{Res}_{y}(p-u, q-v)$ is a polynomial in $\mathbb{C}[x,u,v]$. As (\ref{eq_form_of_p_and_q}) allows us to specialize as noted in Remark \ref{rmk_specialization}, it follows that:
\begin{align}\label{eq-evaluated_resultant}
\operatorname{Res}_{y}(p-c, q-v)=R(x,c,v),
\end{align}
where $R(x,c,v)\in\mathbb{C}[x,v]$ is an irreducible polynomial for every $c\in\mathbb{C}$.

Following the discussion subsequent to Theorem \ref{thm-Adjamagbo_Essen}, we know that to prove the map $(p,q)$ is invertible, it suffices to show that the degree of $x$ in $R(x,u,v)$ is $1$. To achieve this, it is sufficient to prove that $x$ has degree $1$ in $R(x,c,v)$ for every $c\in\mathbb{C}$. To prove this latter fact, we base our analysis -- as in the work of Kraus -- on the theory of algebraic functions within the realm of complex function theory. 

We fix $c\in \mathbb{C}$ and consider the algebraic function $x=f(v)$ determined by the equation $R(x,c,v)=0$; since $R(x,c,v)$ is irreducible, this algebraic function is uniquely determined. Furthermore, $x=f(v)$ is an $n$-valued analytic function on the Riemann sphere -- where $n=\operatorname{deg}_{x}R$ -- except at a finite number of points $\{v_{1}, \dots, v_{\mu}\}\cup\{\infty\}$, where it has, at most, poles or algebraic branch points. Outside these critical values, $x=f(v)$ is represented by $n$ regular function elements, see \cite[Chapter 12]{Hille_62} or \cite[Chapter 8, Section 45]{Markushevich67}. 

Kraus set out to prove that the local condition $J(p,q)=1$ prevents the formation of branch points for $x=f(v)$ in the finite complex plane. However, an algebraic function without branch points in the finite complex plane is necessarily rational; from this, we can conclude that $x$ has degree $1$ in $R(x,c,v)$.

Let's denote $R(x,c,v)=r_{n}(c,v)x^{n}+\cdots+r_{0}(c,v)$. The finite critical points $\{v_{1}, \cdots, v_{\mu}\}$ arise from two sources: (1) $v_{i}$ is a root of the discriminant of the polynomial $R(x,c,v)$ with respect to $x$, provided that $r_{n}(c,v_{i})\ne 0$; or (2) $r_{n}(c,v_{i})=0$. In the first case, although $R(x,c,v_{i})$ has multiple roots, all the roots remain finite and are represented by elements of the form $x=\sum_{l=0}^{\infty}a_{l}(v-v_{i})^{l/\kappa}$, where $1\le\kappa\le n$. In the second case, as $v$ approaches $v_{i}$, one or more roots of $R(x,c,v)$ become infinite and are represented by elements of the form $\frac{1}{x}=\sum_{l=1}^{\infty}a_{l}(v-v_{i})^{l/\kappa}$, where $1\le\kappa\le n$. In both cases, these representations are defined via analytic continuation of regular function elements in a punctured neighborhood of the critical points.

Let $v_{i}$ be a critical point of the first kind, and let $x_{i}$ be one of the finite roots of $R(x,c,v_{i})$, that is, $R(x_{i},c,v_{i})=0$. Then $\operatorname{Res}_{y}(p(x_{i},y)-c, q(x_{i},y)-v_{i})=0$. As $p$ and $q$ are monic in $y$ we can make use of the property of the resultant \cite[Chapter 3, Section 6, Corollary 7]{Cox_Little_Shea_2025} that guarantees the existence of $y_{i}\in\mathbb{C}$ such that $p(x_{i},y_{i})-c=q(x_{i},y_{i})-v_{i}=0$. The curve $p(x,y)-c=0$ is irreducible by hypothesis and smooth due to the condition $J(p,q)=1$. Therefore, we have a local analytic representation of this curve in a neighborhood of $(x_{i},y_{i})$ given by a pair of convergent power series, $x(t)=\sum_{k=0}^{\infty}a_{k}t^{k}$ and $y(t)=\sum_{k=0}^{\infty}b_{k}t^{k}$, satisfying $x(0)=x_{i}$, $y(0)=y_{i}$, $p(x(t),y(t))-c=0$ and $(x'(0),y'(0))\ne (0,0)$.

Define $w(t):=q(x(t),y(t))$ and differentiate it with respect to $t$ as well as the equation $p(x(t),y(t))-c=0$. This yields the following system of linear equations:
\begin{align}\label{eq_linear_system_first_kind}
\frac{\partial p}{\partial x}\frac{\partial x}{\partial t}+\frac{\partial p}{\partial y}\frac{\partial y}{\partial t}=0, \quad \frac{\partial q}{\partial x}\frac{\partial x}{\partial t}+\frac{\partial q}{\partial y}\frac{\partial y}{\partial t}=\frac{\partial w}{\partial t}.
\end{align}
Since $J(p,q)=1$ and $(x'(0),y'(0))\ne (0,0)$, it follows from Cramer's rule applied to (\ref{eq_linear_system_first_kind}) that $w'(0)\ne 0$. Note also that $w(0)=q(x(0),y(0))=q(x_{i},y_{i})=v_{i}$. We can apply the Inverse Function Theorem for analytic functions \cite[Section 9.24]{Nevanlinna-Paatero_69} to write $t$ in terms of $(w-v_{i})$; that is, $t(w)=\sum_{l=1}^{\infty}\alpha_{l}(w-v_{i})^{l}$. Now, we compose this last series with $x(t)=\sum_{k=0}^{\infty}a_{k}t^{k}$  to express $x$ as a convergent power series of $(w-v_{i})$ in a suitable neighborhood of $v_{i}$, that is,
\begin{align}\label{eq_power_series_first_kind}
x(w)=\sum_{l=0}^{\infty}\beta_{l}(w-v_{i})^{l}, \quad x(v_{i})=x_{i}.
\end{align}
On the other hand, it follows from (\ref{eq-evaluated_resultant}) that $R(x(t),c,w(t))=\operatorname{Res}_{y}(p-c, q-q)=0$, this means that (\ref{eq_power_series_first_kind}) solves the equation $R(x,c,w)=0$ in a neighborhood of $v_{i}$. Therefore, $v_{i}$ is a regular point of the algebraic function $x=f(v)$ and it is not a branch point. We have thus proved that there are no critical points of the first kind. 

Observe that the proof above only uses the fact that, at a critical point of the first kind $v_{i}$ the roots $x_{i}$ are finite. By the same argument, we can prove that at a regular point $v$ of the algebraic function $x=f(v)$, the $n$ regular elements can be determined as in (\ref{eq_power_series_first_kind}).

Let $v_{i}$ now be a critical point of the second kind; that is, as $v$ approaches $v_{i}$ one or more roots go to infinity. More precisely, for every $M>0$ there exists a neighborhood of $v_{i}$ such that for every $w_{i}$ in this neighborhood, there exists a root $x_{i}$ of $R(x,c,w_{i})$ such that $||x_{i}||>M$. Since $R(x_{i},c,w_{i})=0$ -- that is, $\operatorname{Res}_{y}(p(x_{i},y)-c, q(x_{i},y)-w_{i})=0$ -- and the polynomials $p$ and $q$ are monic in $y$ we can make use of the property of the resultant \cite[Chapter 3, Section 6, Corollary 7]{Cox_Little_Shea_2025} that ensures the existence of $y_{i}\in\mathbb{C}$ such that $p(x_{i},y_{i})-c=q(x_{i},y_{i})-w_{i}=0$. Putting this all together, we obtain a sequence $\left\{(x_{i},y_{i})\right\}_{i=0}^{\infty}$ such that:
\begin{align}\label{eq_sequence_going_infinity}
||(x_{i},y_{i})||\to \infty \quad \text{and} \quad (p(x_{i},y_{i}),q(x_{i},y_{i}))\to (c,v_{i}).
\end{align}
In view of (\ref{eq_sequence_going_infinity}) we need to take care of the behavior at infinity of the irreducible polynomial $p(x,y)-c=0$, which, by (\ref{eq_form_of_p_and_q}) is monic in $y$. Its Newton-Puiseux expansions at infinity are given by the following fractional power series: $u_{\nu}(x)=x\sum_{k=0}^{\infty}a_{k}\epsilon^{\nu k}x^{-k/m}$, where $m=\operatorname{deg}p$, $\epsilon$ is a primitive $m$-th root of unity, $\nu=0,1,\dots, m-1$, and $\operatorname{gcd}\left(\{k: a_{k}\neq 0\}\cup \{m\}\right)=1$, see Chau \cite[Section 2.1]{Chau99} and Abhyankar \cite[Theorem 5.14]{Abhyankar77}. These series satisfy $p(x,u_{\nu}(x))-c=0$ and provide us with the following meromorphic parametrizations of $p(x,y)-c=0$ for $t$ sufficiently large: $x(t)=t^{m}$, $y(t)=u_{\nu}(t^{m})$, $\nu=0,1,\dots, m-1$.

Define $w_{\nu}(t):=q(x(t),y(t))=q(t^{m},u_{\nu}(t^{m}))$, $\nu=0,1,\dots, m-1$. Observe that each $w_{\nu}(t)$ has at most finitely many positive powers of $t$. Since we already know the expansion at infinity of $p(x,y)-c=0$ and $p(t^{m},u_{\nu}(t^{m}))-c=0$, the existence of the sequence $\left\{(x_{i},y_{i})\right\}_{i=0}^{\infty}$ satisfying (\ref{eq_sequence_going_infinity}) implies -- possibly after passing to a subsequence -- that there exists at least one among the $w_{\nu}(t)$, say $w_{\alpha}(t)$, such that $w_{\alpha}(t)=\sum_{k=0}^{\infty}b_{k}t^{-k}$ and $b_{0}=v_{i}$.

Make the change of variables $s=1/t$ to bring the neighborhood of infinity to the neighborhood of zero. Differentiating with respect to $s$ the equations $p(x(s),y(s))-c=0$ and $w_{\alpha}(s)=q(x(s),y(s))$, we obtain the following system of linear equations: 
\begin{align}\label{eq_linear_system_second_kind}
\frac{\partial p}{\partial x}\frac{\partial x}{\partial s}+\frac{\partial p}{\partial y}\frac{\partial y}{\partial s}=0,\quad \frac{\partial q}{\partial x}\frac{\partial x}{\partial s}+\frac{\partial q}{\partial y}\frac{\partial y}{\partial s}=\frac{\partial w_{\alpha}}{\partial s}.
\end{align} 
\noindent
{\itshape Kraus erroneously claimed \cite[p. 822]{Kraus1884} that -- since $x'(0)\ne 0$ (because $x(s)=\frac{1}{s^{m}}$ and $x'(s)=\frac{-m}{s^{m+1}}$) and $J(p,q)=1$ -- we can conclude that $w_{\alpha}'(0)\ne 0$ by solving the linear system above. However, it should be noted that the derivative of $x(s)$ at $s=0$ is not defined.}

If we were able to prove that $w_{\alpha}'(0)\ne 0$, we could proceed as in the case of a critical point of the first kind. We could apply the Inverse Function Theorem for analytic functions, and since $w_{\alpha}(0)=b_{0}=v_{i}$ we could write $s$ in terms of $(w_{\alpha}-v_{i})$; that is, $s(w_{\alpha})=\sum_{l=1}^{\infty}d_{l}(w_{\alpha}-v_{i})^{l}$. Now, we would compose this last series with $x(t)=t^{m}=\frac{1}{s^{m}}$  to express $\frac{1}{x}$ as a convergent power series of $(w_{\alpha}-v_{i})$ in some suitable neighborhood of $v_{i}$, that is,
\begin{align}\label{eq_power_series_second_kind}
	\frac{1}{x(w_{\alpha})}=\sum_{l=1}^{\infty}\gamma_{l}(w_{\alpha}-v_{i})^{l}.
\end{align}
On the other hand, it would follow from (\ref{eq-evaluated_resultant}) that $R(x(t),c,w_{\alpha}(t))=\operatorname{Res}_{y}(p-c, q-q)=0$. This would mean that (\ref{eq_power_series_second_kind}) solves the equation $R(x,c,w)=0$ in a neighborhood of $v_{i}$. Therefore, $v_{i}$ would be a regular point of the algebraic function $x=f(v)$ and not a branch point. We would have proved that there are no critical points of the second kind either. Unfortunately, we don't know if $w_{\alpha}'(0)\ne 0$.

Note that the existence of critical points of the second kind -- that is, roots of $r_{n}(c,v)=0$ -- is precisely the source of the non-properness of the map $(p,q)$. This is reflected in (\ref{eq_sequence_going_infinity}); we could have a sequence in the domain going to infinity while its images converge to a finite point.  In fact, the set of points at which $(p,q)$ is not proper is exactly the set $\{(u,v)\in\mathbb{C}^{2} \mid r_{n}(u,v)=0\}$, see Jelonek \cite[Theorem 2.1]{Jelonek01}.

Finally, we would like to mention a theorem of Chadzyński and Krasi{\'n}ski \cite{Chadzynski_Krasinski_93}. In their proof, they use precisely the system of linear equations (\ref{eq_linear_system_second_kind}) to show that the map $(p,q)$ is invertible, provided that its restriction to the set $\{(z_1,z_2)\in\mathbb{C}^{2} \mid p(z_1,z_2)=0\}$ is proper.

\bibliographystyle{abbrv}
\bibliography{bibliografia}

@misc{Weierstrass_1902,
	author = {Weierstra{\ss}, K.},
	title = {Mathematische {Werke}. {Herausgegeben} unter {Mitwirkung} einer von der k{\"o}niglich preu{{\ss}}ischen {Akademie} der {Wissenschaften} eingesetzten {Kommission}. {Vierter} {Band}. {Vorlesungen} {\"u}ber die {Theorie} der {Abelschen} {Transzendenten}. {Bearbeitet} von \emph{{G}. {Hettner}} und \emph{{J}. {Knoblauch}}.},
	year = {1902},
	language = {German},
	howpublished = {Berlin: {Mayer} {{\(\&\)}} {M{\"u}ller}. {XIV} u. 631 {S}. {{\(4^{\circ}\)}} (1902).},
	zbMATH = {2658832},
	JFM = {33.0031.01}
}

@incollection {Moh98,
	AUTHOR = {Moh, T. T.},
	TITLE = {Jacobian conjecture},
	BOOKTITLE = {Algebra and geometry ({T}aipei, 1995)},
	SERIES = {Lect. Algebra Geom.},
	VOLUME = {2},
	PAGES = {103--116},
	PUBLISHER = {Int. Press, Cambridge, MA},
	YEAR = {1998},
	ISBN = {1-57146-058-6},
	MRCLASS = {14R15},
	MRNUMBER = {1697951},
	MRREVIEWER = {L.\ Andrew\ Campbell},
}

@article {Jelonek01,
	AUTHOR = {Jelonek, Zbigniew},
	TITLE = {Note about the set {{\(S_f\)}} for a polynomial mapping {{\(f: {\mathbb{C}}^2{{\to}} {\mathbb{C}}^2\)}}},
	JOURNAL = {Bull. Polish Acad. Sci. Math.},
	FJOURNAL = {Polish Academy of Sciences. Bulletin. Mathematics},
	VOLUME = {49},
	YEAR = {2001},
	NUMBER = {1},
	PAGES = {67--72},
	ISSN = {0239-7269},
	MRCLASS = {14R10},
	MRNUMBER = {1824157},
	MRREVIEWER = {James\ K.\ Deveney},
}

@article {Chadzynski_Krasinski_93,
	AUTHOR = {Chadzyński, Jacek and Krasi{\'n}ski, Tadeusz},
	TITLE = {Properness and the {Jacobian} conjecture in {{\(\mathbb{C}^2\)}}},
	JOURNAL = {Bull. Soc. Sci. Lett. \L \'{o}d\'{z} S\'{e}r. Rech.
	D\'{e}form.},
	FJOURNAL = {Bulletin de la Soci\'{e}t\'{e} des Sciences et des Lettres de
	\L \'{o}d\'{z}. S\'{e}rie: Recherches sur les
	D\'{e}formations},
	VOLUME = {14},
	YEAR = {1992/93},
	NUMBER = {131-140},
	PAGES = {13--19},
	ISSN = {0459-6854,2450-9329},
	MRCLASS = {14E07 (32H02 32H04 32H35)},
	MRNUMBER = {1255660},
	MRREVIEWER = {L.\ Andrew\ Campbell},
}

@article {Chadzynski_Krasinski_92,
	AUTHOR = {Chadzyński, Jacek and Krasi{\'n}ski, Tadeusz},
	TITLE = {On the {{\L}}ojasiewicz exponent at infinity for polynomial mappings of {{\(\mathbb{C}^2\)}} into {{\(\mathbb{C}^2\)}} and components of polynomial automorphisms of {{\(\mathbb{C}^2\)}}},
	JOURNAL = {Ann. Polon. Math.},
	FJOURNAL = {Annales Polonici Mathematici},
	VOLUME = {57},
	YEAR = {1992},
	NUMBER = {3},
	PAGES = {291--302},
	ISSN = {0066-2216,1730-6272},
	MRCLASS = {14E05 (14E09)},
	MRNUMBER = {1201856},
	MRREVIEWER = {Zbigniew\ Jelonek},
	DOI = {10.4064/ap-57-3-291-302},
	URL = {https://doi.org/10.4064/ap-57-3-291-302},
}

@article {Bass_Connell_Wright_82,
	AUTHOR = {Bass, Hyman and Connell, Edwin H. and Wright, David},
	TITLE = {The {J}acobian conjecture: reduction of degree and formal
	expansion of the inverse},
	JOURNAL = {Bull. Amer. Math. Soc. (N.S.)},
	FJOURNAL = {American Mathematical Society. Bulletin. New Series},
	VOLUME = {7},
	YEAR = {1982},
	NUMBER = {2},
	PAGES = {287--330},
	ISSN = {0273-0979,1088-9485},
	MRCLASS = {14H20 (13B25)},
	MRNUMBER = {663785},
	MRREVIEWER = {Stephen\ McAdam},
	DOI = {10.1090/S0273-0979-1982-15032-7},
	URL = {https://doi.org/10.1090/S0273-0979-1982-15032-7},
}

@article{Jacobi1841,
	author = {Jacobi, C. G. J.},
	title = {De {Determinantibus} functionalibus.},
	fjournal = {Journal f{\"u}r die Reine und Angewandte Mathematik},
	journal = {J. Reine Angew. Math.},
	issn = {0075-4102},
	volume = {22},
	pages = {319--359},
	year = {1841},
	language = {Latin},
	doi = {10.1515/crll.1841.22.319},
	url = {https://eudml.org/doc/147138},
	zbMATH = {2751099},
	ERAM = {022.0686cj}
}

@book{Jacobi1896,
	author = {Jacobi, C. G. J.},
	title = {Ueber die {Functionaldeterminanten} (1841). {Hersg}. von {P}. {St{\"a}ckel}.},
	fseries = {Ostwalds Klassiker der Exakten Wissenschaften},
	series = {Ostwalds Klass. Exakten Wiss.},
	volume = {78},
	year = {1896},
	publisher = {Harri Deutsch, Frankfurt am Main},
	language = {German},
	zbMATH = {2675614},
	JFM = {27.0107.03}
}

@misc{Kraus1883,
	author = {Kraus, L.},
	title = {Ueber rational umkehrbare {Substitutionen}.},
	year = {1883},
	language = {German},
	howpublished = {Prag. {Ber}. 187-196},
	zbMATH = {2701493},
	JFM = {16.0065.02}
}

@misc{Kraus1882,
	author = {Kraus, L.},
	title = {Ueber rational umkehrbare {Substitutionen}.},
	year = {1882},
	language = {German},
	howpublished = {Prag. {Ber}. 338-352},
	zbMATH = {2703371},
	JFM = {15.0113.02}
}

@article{Weyr1886,
	author    = {Eduard Weyr},
	title     = {{Život a působení dra Ludvíka Krause} ({C}zech) [{L}ife and activities of {D}r. {L}udvík {K}raus]},
	journal   = {{Časopis pro pěstování mathematiky a fysiky}},
	volume    = {15},
	year      = {1886},
	number    = {2},
	pages     = {[49a]--52},
	language  = {Czech}
}

@misc{Muir1923,
	author = {Muir, T.},
	title = {The theory of determinants in the historical order of development. {Vol}. 4: {The} period 1880 to 1900.},
	year = {1923},
	language = {English},
	howpublished = {London: {Macmillan}, 31 u. 508 {S}. {{\(8^\circ\)}}},
	zbMATH = {2597826},
	JFM = {49.0078.01}
}

@article{Kraus1884,
	author = {Kraus, L.},
	title = {Ueber {Functionaldeterminanten}.},
	fjournal = {Akademie der Wissenschaften in Wien, Mathematisch-Naturwissenschaftliche Klasse, Sitzungsberichte, Abteilung IIa},
	journal = {Wien. Ber.},
	volume = {90},
	pages = {813--826},
	year = {1884},
	language = {German},
	zbMATH = {2701862},
	JFM = {16.0232.03}
}

@book {Abhyankar77,
	AUTHOR = {Abhyankar, S. S.},
	TITLE = {Lectures on expansion techniques in algebraic geometry},
	SERIES = {Tata Institute of Fundamental Research Lectures on Mathematics
	and Physics},
	VOLUME = {57},
	NOTE = {Notes by Balwant Singh},
	PUBLISHER = {Tata Institute of Fundamental Research, Bombay},
	YEAR = {1977},
	PAGES = {iv+168},
	MRCLASS = {14H20},
	MRNUMBER = {542446},
	MRREVIEWER = {Jos\'e\ L.\ Vicente},
}

@book {Markushevich67,
	AUTHOR = {Markushevich, A. I.},
	TITLE = {Theory of functions of a complex variable. {V}ol. {III}},
	NOTE = {Revised English edition, translated and edited by Richard A.
	Silverman},
	PUBLISHER = {Prentice-Hall, Inc., Englewood Cliffs, NJ},
	YEAR = {1967},
	PAGES = {xi+360},
	MRCLASS = {30.00},
	MRNUMBER = {215964},
	MRREVIEWER = {R. P. Boas},
}

@article {Ploski94,
	AUTHOR = {P{\l}oski, Arkadiusz},
	TITLE = {A note on the {{\L}}ojasiewicz exponent at infinity},
	JOURNAL = {Bull. Soc. Sci. Lett. \L \'{o}d\'{z} S\'{e}r. Rech.
	D\'{e}form.},
	FJOURNAL = {Bulletin de la Soci\'{e}t\'{e} des Sciences et des Lettres de
	\L \'{o}d\'{z}. S\'{e}rie: Recherches sur les
	D\'{e}formations},
	VOLUME = {17},
	YEAR = {1994},
	PAGES = {11--15},
	ISSN = {0459-6854,2450-9329},
	MRCLASS = {32H35 (14E05)},
	MRNUMBER = {1350859},
}

@article {Sakkalis93,
	AUTHOR = {Sakkalis, Takis},
	TITLE = {On relations between {J}acobians and resultants of polynomials
	in two variables},
	JOURNAL = {Bull. Austral. Math. Soc.},
	FJOURNAL = {Bulletin of the Australian Mathematical Society},
	VOLUME = {47},
	YEAR = {1993},
	NUMBER = {3},
	PAGES = {473--481},
	ISSN = {0004-9727},
	MRCLASS = {12D10 (12E05 13B10 14E07)},
	MRNUMBER = {1220322},
	MRREVIEWER = {L.\ Andrew\ Campbell},
	DOI = {10.1017/S0004972700015306},
	URL = {https://doi.org/10.1017/S0004972700015306},
}

@book {Nevanlinna-Paatero_69,
	AUTHOR = {Nevanlinna, Rolf and Paatero, V.},
	TITLE = {Introduction to complex analysis},
	NOTE = {Translated from the German by T. K\"ovari and G. S. Goodman},
	PUBLISHER = {Addison-Wesley Publishing Co., Reading, Mass.-London-Don
	Mills, Ont.},
	YEAR = {1969},
	PAGES = {ix+348},
	MRCLASS = {30.00},
	MRNUMBER = {239056},
}

@book {Hille_62,
	AUTHOR = {Hille, Einar},
	TITLE = {Analytic function theory. {V}ol. {II}},
	SERIES = {Introductions to Higher Mathematics},
	PUBLISHER = {Ginn and Company, Boston, Mass.-New York-Toronto},
	YEAR = {1962},
	PAGES = {xii+496},
	MRCLASS = {30.00},
	MRNUMBER = {201608},
	MRREVIEWER = {G.\ R.\ MacLane},
}

@article {Adjamagbo_Essen90,
	AUTHOR = {Adjamagbo, Kossivi and van den Essen, Arno},
	TITLE = {A resultant criterion and formula for the inversion of a
	polynomial map in two variables},
	JOURNAL = {J. Pure Appl. Algebra},
	FJOURNAL = {Journal of Pure and Applied Algebra},
	VOLUME = {64},
	YEAR = {1990},
	NUMBER = {1},
	PAGES = {1--6},
	ISSN = {0022-4049},
	MRCLASS = {14E07 (13B10 14E20)},
	MRNUMBER = {1055017},
	MRREVIEWER = {L. Andrew Campbell},
	DOI = {10.1016/0022-4049(90)90002-Y},
	URL = {https://doi.org/10.1016/0022-4049(90)90002-Y},
}

@article {Cassou-Nogues_Miyanishi09,
	AUTHOR = {Cassou-Nogu\`es, Pierrette and Miyanishi, Masayoshi},
	TITLE = {Smoothness of the images of members of a linear system under
	an endomorphism of the affine plane},
	JOURNAL = {J. Pure Appl. Algebra},
	FJOURNAL = {Journal of Pure and Applied Algebra},
	VOLUME = {213},
	YEAR = {2009},
	NUMBER = {5},
	PAGES = {711--723},
	ISSN = {0022-4049},
	MRCLASS = {14R10 (14R15)},
	MRNUMBER = {2494363},
	MRREVIEWER = {Yuriy Bodnarchuk},
	DOI = {10.1016/j.jpaa.2008.09.002},
	URL = {https://doi.org/10.1016/j.jpaa.2008.09.002},
}

@book{Cox_Little_Shea_2025,
	author = {Cox, David A. and Little, John and O'Shea, Donald},
	title = {Ideals, varieties, and algorithms. {An} introduction to computational algebraic geometry and commutative algebra},
	edition = {5th edition},
	isbn = {978-3-031-91840-7; 978-3-031-91843-8; 978-3-031-91841-4},
	year = {2025},
	publisher = {Cham: Springer},
	language = {English},
	zbMATH = {8074064}
}

@book {Abhyankar71,
	AUTHOR = {Abhyankar, Shreeram S.},
	TITLE = {Algebraic space curves},
	SERIES = {S\'eminaire de Math\'ematiques Sup\'erieures [Seminar on
	Higher Mathematics]},
	VOLUME = {No. 43 (\'Et\'e 1970)},
	PUBLISHER = {Les Presses de l'Universit\'e{} de Montr\'eal, Montreal, QC},
	YEAR = {1971},
	PAGES = {114},
	MRCLASS = {14H45 (14N10)},
	MRNUMBER = {399109},
	MRREVIEWER = {P.\ Schenzel},
}

@article {Kaliman93,
    AUTHOR = {Kaliman, Shulim},
     TITLE = {On the {J}acobian conjecture},
   JOURNAL = {Proc. Amer. Math. Soc.},
  FJOURNAL = {Proceedings of the American Mathematical Society},
    VOLUME = {117},
      YEAR = {1993},
    NUMBER = {1},
     PAGES = {45--51},
      ISSN = {0002-9939,1088-6826},
   MRCLASS = {14E07 (13B10 14E20)},
  MRNUMBER = {1106179},
MRREVIEWER = {L.\ Andrew\ Campbell},
       DOI = {10.2307/2159696},
       URL = {https://doi.org/10.2307/2159696},
}

@article {Keller39,
	AUTHOR = {Keller, Ott-Heinrich},
	TITLE = {Ganze {C}remona-{T}ransformationen},
	JOURNAL = {Monatsh. Math. Phys.},
	FJOURNAL = {Monatshefte f\"ur Mathematik und Physik},
	VOLUME = {47},
	YEAR = {1939},
	NUMBER = {1},
	PAGES = {299--306},
	ISSN = {1812-8076},
	MRCLASS = {99-04},
	MRNUMBER = {1550818},
	DOI = {10.1007/BF01695502},
	URL = {https://doi.org/10.1007/BF01695502},
}

@article {Formanek94,
	AUTHOR = {Formanek, Edward},
	TITLE = {Observations about the {J}acobian conjecture},
	JOURNAL = {Houston J. Math.},
	FJOURNAL = {Houston Journal of Mathematics},
	VOLUME = {20},
	YEAR = {1994},
	NUMBER = {3},
	PAGES = {369--380},
	ISSN = {0362-1588},
	MRCLASS = {14E09 (13F30 14E05)},
	MRNUMBER = {1287981},
	MRREVIEWER = {Ludwik M. Dru\.{z}kowski},
}

@article {McKay_Wang86,
	AUTHOR = {McKay, James H. and Wang, Stuart Sui Sheng},
	TITLE = {An inversion formula for two polynomials in two variables},
	JOURNAL = {J. Pure Appl. Algebra},
	FJOURNAL = {Journal of Pure and Applied Algebra},
	VOLUME = {40},
	YEAR = {1986},
	NUMBER = {3},
	PAGES = {245--257},
	ISSN = {0022-4049},
	MRCLASS = {12E05},
	MRNUMBER = {836651},
	MRREVIEWER = {Ira Gessel},
	DOI = {10.1016/0022-4049(86)90044-7},
	URL = {https://doi.org/10.1016/0022-4049(86)90044-7},
}

@article {Gutierrez_Rubio_Yu02,
	AUTHOR = {Gutierrez, Jaime and Rubio, Rosario and Yu, Jie-Tai},
	TITLE = {{$D$}-resultant for rational functions},
	JOURNAL = {Proc. Amer. Math. Soc.},
	FJOURNAL = {Proceedings of the American Mathematical Society},
	VOLUME = {130},
	YEAR = {2002},
	NUMBER = {8},
	PAGES = {2237--2246},
	ISSN = {0002-9939},
	MRCLASS = {13P99 (68W30)},
	MRNUMBER = {1896403},
	MRREVIEWER = {Ming Zhang},
	DOI = {10.1090/S0002-9939-02-06331-1},
	URL = {https://doi.org/10.1090/S0002-9939-02-06331-1},
}

@article {Hejmej18,
	AUTHOR = {Hejmej, Beata},
	TITLE = {A note about irreducibility of a resultant},
	JOURNAL = {Bull. Soc. Sci. Lett. \L \'{o}d\'{z} S\'{e}r. Rech. D\'{e}form.},
	FJOURNAL = {Bulletin de la Soci\'{e}t\'{e} des Sciences et des Lettres de \L \'{o}d\'{z}.
	S\'{e}rie: Recherches sur les D\'{e}formations},
	VOLUME = {68},
	YEAR = {2018},
	NUMBER = {1},
	PAGES = {27--31},
	ISSN = {0459-6854},
	MRCLASS = {13P15 (12F10)},
	MRNUMBER = {4134733},
	DOI = {10.1093/imrn/rnu106},
	URL = {https://doi.org/10.1093/imrn/rnu106},
}

@article {Garcia_Barroso-Gwozdziewic22,
	AUTHOR = {Garc\'{\i}a Barroso, Evelia R. and Gwo\'{z}dziewicz, Janusz},
	TITLE = {Higher order polars of quasi-ordinary singularities},
	JOURNAL = {Int. Math. Res. Not. IMRN},
	FJOURNAL = {International Mathematics Research Notices. IMRN},
	YEAR = {2022},
	NUMBER = {2},
	PAGES = {1045--1080},
	ISSN = {1073-7928},
	MRCLASS = {14H20 (14M25 32B05)},
	MRNUMBER = {4368879},
	MRREVIEWER = {Angel Granja},
	DOI = {10.1093/imrn/rnaa106},
	URL = {https://doi.org/10.1093/imrn/rnaa106},
}

@article {Chau99,
	AUTHOR = {Nguyen Van Chau},
	TITLE = {Non-zero constant {J}acobian polynomial maps of {${\bf C}^2$}},
	JOURNAL = {Ann. Polon. Math.},
	FJOURNAL = {Annales Polonici Mathematici},
	VOLUME = {71},
	YEAR = {1999},
	NUMBER = {3},
	PAGES = {287--310},
	ISSN = {0066-2216,1730-6272},
	MRCLASS = {14R15},
	MRNUMBER = {1704304},
	MRREVIEWER = {Ludwik\ M.\ Dru\.zkowski},
	DOI = {10.4064/ap-71-3-287-310},
	URL = {https://doi.org/10.4064/ap-71-3-287-310},
}

\end{document}